\documentclass[12pt]{article}
\usepackage{amssymb}
\newtheorem{qwe}{Statement}
\newtheorem{Theorem}{Theorem}[section]
\newtheorem{Corollary}{Corollary}[section]
\newtheorem{Lemma}{Lemma}[section]
\newtheorem{Statement}{Statement}
\newtheorem{Definition}{Definition}
\newtheorem{Sta}{Statement}[section]
\newtheorem{Def}{Definition}[section]
\begin{document}

\title{On Special monodromy groups and Riemann-Hilbert problem for Riemann equation}

\author{Vladimir Poberezhny \thanks{poberezh@itep.ru}
\date{} \\
{\small {\it Steklov Mathemayical Institute}}
}
\maketitle


Riemann equation is a second order Fuchsian differential equation having
three singular points ${\cal D}=\{a_1,a_2,a_3\}$ on a Riemann sphere ${\mathbb
{CP}}^1$. As any linear differential equation it produce some monodromy
representation. Taking a gerbe of fundamental matrix $Y(z)$ in a neighborhood
 of point $z_0$ , and continuing it along a loop $\gamma$ we obtain
another fundamental matrix $Y^*(z)$. Beeing also a gerbe of fundamental set
of solutions $Y^*(z)$ is connected with $Y(z)$ via multiplication on some
constant matrix $Y(z)=Y^*(z)G_{\gamma}$. It is easy to see that $G_{\gamma}$
is defined not by the loop $\gamma$, but by its homotopy class $[\gamma]$.
This homomorphism $\chi:\pi_1({\mathbb CP}^1\setminus{\cal D})\to
GL(2,{\mathbb C})$ is called monodromy representation and its image is
the monodromy of equation.

The problem of constructing an equation with a given monodromy was posed for
the first time by Riemann. The question of existence of Fuchsian equation
with a given monodromy and prescribed singularities was included by Hilbert,
in his list of mathematical problems. Now it is known as 21st Hilbert problem
or Riemann-Hilbert Problem.

In the present work we explore monodromy of Riemann equation and solve
Riemann-Hilbert problem for it. It is known that  a Fuchsian equation of order
$p$ having $n$ singular points depends on $\frac{p^2}{2}(n-2)+\frac{pn}{2}$
parameters. The corresponding monodromy representation is described by
$p^2(n-2)+1$ parameters. One can see  that in general if $n>3$ and $p>2$
then it is impossible to construct an equation with a given monodromy,
since there exist more monodromy
representations than fuchsian equations. So, the case of
Riemann equation $(p=2,n=3)$ is in some sense distinguished.

In the first part of the present work we use the Levelt's valuation theory to
describe all monodromy representations that can be realized by Riemann
equation. In the second part we show that if the monodromy of Riemann
equation lies in $SL(2,\mathbb{C})$, then such a monodromy  can also be
realized by a more special Riemann-Sturm-Liouville equation. The third part is
devoted to $SL(2,\mathbb{Z})$ monodromy. We construct the criterion for
hypergeometric equation to have monodromy in $SL(2,\mathbb{Z})$. The results
of first part were obtained earlier in \cite{yap,kim} using Fuchs relation and
algebraic methods of representation theory.

This paper consists of a collection of student works done by author during his
studies at Moscow State University under tutoring of \fbox{Andrey
Bolibruch}. His support, attention and useful discussions were unvaluable.

The question of the first part,
that is Riemann-Hilbert problem for Riemann equation, was also given to
S.~Malek, another student of Andrey Bolibruch. He independently succeded in its
solution as well.

The work is supported by grants RFBR-04-01-00642, NS-457.2003.1, PICS-2094 and
by the program "Mathematical methods in non-linear dynamics".

\section*{Preliminaries}
In this section we mention some facts regarding Fuchsian equation, that
will be of use later. We can without loss of generality assume zero as
one of singular points of equation to give some results from local theory.
\begin{Definition} The equation
$$
y^{(k)}+b_1(z)y^{(k-1)}+...+b_k(z)y=0
$$
is Fuchsian at zero, if each $b_j(z)$ have there a pole of order at most $j$.
In other words, $b^jx^j$ is holomorphic in some neighborhood of zero.
\end{Definition}
\begin{Definition} The system
$$ \dot y=B(z)y\qquad
\mbox{$y\in {\mathbb C}^n$,\,\,$B(z)\in M_n({\mathbb C})$ }
$$
is Fuchsian at its singular point zero if  $B(z)$ has a simple pole at zero.
\end{Definition}
\begin{Definition} The system (equation) is fuchsian, if  it is fuchsian
at all its singular points.
\end{Definition}
\begin{Definition} Zero is a regular singular point of an
equation (system) if any component of any its solution is of at most polynomial
growth when $z$ tends to $0$, in any sector of finite angle having vertex at
zero.
\end{Definition}
There are some relations between notions of fuchsian and regular equations
and systems.
\begin{Statement} (\cite{b7,zog})
The equation is Fuchsian iff it is regular.
\end{Statement}
\begin{Statement} (\cite{zog,for})
If the system is Fuchsian, then it is regular.
\end{Statement}

For further investigations of asymptotics of an equation at its singular
point we will use the following notion of valuation of the function due
to Levelt.

\begin{Definition} Levelt's norm of function $f(z)$ at zero is
defined as follows:
$$
\varphi_0(f)=\sup\left(\lambda\in {\mathbb Z}:\frac{f(z)}{z^\lambda}\to
0\quad\mbox{when $z\to 0$ in a sector of span less than $2\pi$ }\right)
$$
The valuation of a vector or matrix is defined as the minimum of
the valuations of their components.
\end{Definition}
In particular for $f(z)=\frac{c_{-k}}{z^k}+...+c_0+...+c_jz^j+...$ we have
 $\varphi_0(f)=-k$.

Described valuation  has the following properties:
\begin{enumerate}
\item $\varphi_0(c\cdot f)=\varphi_0(f)$ for any $c\in{\mathbb C}$.
\item $\varphi_0(f\cdot\ln^\alpha z)=\varphi_0(f)$ for any $\alpha\in{\mathbb C}$.
\item $\varphi_0(f_1+f_2)\ge\min\left(\varphi_0(f_1),\varphi_0(f_2)\right)$,
if $\varphi_0(f_1)\ne\varphi_0(f_2)$ it is a strong equality.
 \item $\varphi_0(g^*f)=\varphi_0(f)$ where $g^*$- is a monodromy operator
corresponding to the counter-clockwise passing around zero.
\end{enumerate}

The Levelt's valuation define subspaces ${\cal X}^j$ of the solution space $\cal X$.
$$
{\cal X}^j=\{ f\in{\cal X}: \varphi_0(f)\ge j\}
$$
We have ${\cal X}^i\subseteq{\cal X}^k$ for $i\ge k$.
As ${\cal X}$ is finite-dimensional we see that valuation can take only a
finite number of values $-\infty<\varphi_0^1<...<\varphi_0^k<\infty$.
That defines the Levelt's filtration:
$$
0\subset{\cal X}^{\varphi_k}....\subset{\cal X}^{\varphi_1}={\cal X}
$$
One can note, that ${\cal X}^{\varphi_{j}}$ are invariant under monodromy
operator action. It follows from the Levelt's valuation properties.

\begin{Definition} Consider a basis $(e_1,...,e_p)$ of the solution space ${\cal X}$.
This basis is weakly associated at singular point zero if valuation $\varphi_0(\cdot)$
takes on it all its values with multiplicities.
\end{Definition}
\begin{Definition} The weakly associated basis $(e_1,...,e_p)$ is associated at zero
if matrix of monodromy operator $G_{0}$ is upper-triangular in that basis and
 $\varphi_{0}(e_k)\ge\varphi_{0}(e_{k+1})$.
\end{Definition}
\begin{Definition} The basis fo solution space ${\cal X}$ constructed as the union
of associated basises of eigenspaces of monodromy operator is called Levelt's basis.
 \end{Definition}
\begin{Statement}(\cite{Bol})
Levelt's basis is also weakly associated.
\end{Statement}
\begin{Definition} Take $G$ to be a matrix of monodromy operator.
Then its normalized logarithm is the following matrix:
$$
E=\frac{1}{2\pi i}\ln G\qquad \mbox{$\Re\rho_i\in[0;1)$ }
$$
 where
$\rho_i$- are the the eigenvalues of $E$, аnd $\Re\rho_i$ are their real parts.
\end{Definition}
One can easily check that after analytical continuation around zero we have
 $z^{-E_0}\to G_0^{-1}z^{-E_0}$. So the following decomposition of the fundamental
matrix in a zero's neighborhood holds
$$
Y=M(z)z^{E_0}
$$
Here $M(z)$ is a meromorphic matrix function, and $E_0$
is the corresponding normalized logarithm.
\begin{Statement}\label{nevyr}(\cite{lev})
For fundamental matrix $Y(z)$ constructed from the basis of ${\cal X}$
weakly associated at regular singular point 0  the following decomposition holds.
$$
Y(z)=U_0(z)z^{A_0}z^{E_0}
$$
where $E_0$ is the normalized logarithm of $G_0$,
$A_0$ is an integer diagonal matrix of valuations of the basis
$A_0=diag (\varphi_0^1,...\varphi_0^p)$, and $U_0(z)$ is holomorphic at 0.
 The system is Fuchsian at 0 iff $\det U_0(0)\ne 0$.
\end{Statement}
\begin{Definition}
The numbers
$\beta_i^j=\varphi_{i}^j+\rho_i^j$
being on the diagonal of matrix $A_{a_i}+E_{a_i}$, constructed in
the basis associated at $a_i$ are the exponents of the space ${\cal X}$
at the point $a_i$.
\end{Definition}
 \begin{Statement}(\cite{Bol})
The sum of exponents of a system having only regular singular points is a
non-positive integer
$$
\sum_{\cal D}\sum_{j=1}^{j=p}\beta_i^j=k\le0
$$
The system is Fuchsian iff
$$
\sum_{\cal D}\sum_{j=1}^{j=p}\beta_i^j=0
$$
 \end{Statement}
\begin{Statement}\label{fuks}(\cite{gol})
For a scalar fuchsian equation of order $p$ having $n$ singular points the
following equality holds
$$
\sum_{\cal D}\sum_{j=1}^{j=p}\beta_i^j=\frac {(n-2)p(p-1)}{2}
$$
\end{Statement}
The problem of constructing Fuchsian system having a prescribed
monodromy is known as the Riemann-Hilbert problem. This work is based on the
following result of Krylov (\cite{kr}).
\begin{Statement} \label{kr} In the class of systems of second rank with a three
singular points the Rieman-Hilbert problem can be solved  positively.
In other words, there always exists a fuchsian system of rank 2 with three given
singular points and any prescribed monodromy.
\end{Statement}
That result was generalized later by Dekkers to an arbitrary number of
singularities.(\cite{dek}).

\section{$GL(2,{\mathbb C})$ monodromy}
\subsection{Irreducible representations}
Let $\chi$ be an irreducible $GL(2,{\mathbb C})$-representation
and take $Y$ --  a fundamental matrix realizing  $\chi$ (such a matrix
always exists, according to \ref{kr}). We can think of $Y$ as a matrix
representing  Levelt's basis at $a_i$ for some $i$.
$$
 Y(z)=U_i(z)z^{A_i}z^{E_i}\qquad\mbox{$\det U_i(a_i)\ne 0$}
$$
Consider
$$
\widetilde Y=
  \left(
    \begin{array}{cc}
      0 & 1 \\
      1 & 0
    \end{array}
  \right)
U_i^{-1}(a_i)Y
$$
Evidently, $\widetilde Y$ is a fundamental matrix of a fuchsian system with
the same monodromy and singularities. We will denote it by $Y$ too.
\begin {Lemma}{ If $\chi$ is irreducible, then the elements of the first row of
the matrix $Y(z)$, i.e. functions $y_1(z),y_2(z)$, are linearly independent.}
\end{Lemma}
{\bf Proof:} Suppose $y_1(z),y_2(z)$
are dependent, this implies
\begin{equation}
y_2(z)=ky_1(z)\label{zav}
\end{equation}
for some $k$. That leads to some relations for matrices $G_{a_i}$.
Since from (\ref{zav}) we have $g^*_i(y_2)=kg^*_i(y_1)$ for all  $i$,
where $g^*_i$ is a monodromy operator corresponding to loop around $a_i$,
then, for all $G_{a_i}$ of the type
$$
G_{a_i}=
  \left(
    \begin{array}{cc}
      \alpha_i & \beta_i\\
      \gamma_i &\delta_i
    \end{array}
  \right)
$$
we have
$$
k\alpha_i+k^2\gamma_i=\beta_i+k\delta_i
$$
Taking now
$$
S=
  \left(
    \begin{array}{cc}
      k&0\\
      1&1
    \end{array}
  \right),
S^{-1}=
  \left(
    \begin{array}{cc}
      k^{-1}&0\\
      -k^{-1}&1
    \end{array}
  \right)
$$
and turning to the basis $\tilde Y=YS$, we get
$$
\widetilde G_{a_i}=S^{-1}G_{a_i}S=
  \left(
    \begin{array}{cc}
      \alpha_i-k^{-1}\beta_i&k^{-1}\\
      0&\delta_i-k^{-1}\beta_i
    \end{array}
  \right)
$$
for any  $i$. That contradicts to irreducibility of $\chi$. $\Box$

So, let $Y$ be a matrix representing a Levelt's basis, such that
 $$
Y=
  \left(
    \begin{array}{cc}
      y_1 & y_2 \\
      y_3 & y_4
    \end{array}
  \right)\qquad
U_i(a_i)=
  \left(
    \begin{array}{cc}
      0 & 1 \\
      1 & 0
    \end{array}
  \right)
$$
We take its first line  $(y_1,y_2)$ as the basis of some scalar regular
equation, and we will show that this equation is exactly the Riemann
equation. From the form of $U_i(a_i)$ it follows that
$\varphi_{a_i}(y_1)\ge\varphi_{a_i}^1$.
The exponents of the scalar equation at points $\cal D$ are evidently
estimated by the exponents of the system:
$$
\tilde \beta_i^1\ge\beta_i^1+1,\,\,\tilde\beta_i^2=\beta_i^2
$$
$$
\Re\tilde\beta_j^k\ge\Re\beta_j^k,\,\, \Im\tilde\beta_j^k=\Im\beta_j^k,\,\,j\ne i
$$
\begin{Lemma}\label{postrur}\,\, Let $(\tilde y_1,\tilde y_2)$ be linearly independent
 and holomorphic on $({\overline{\mathbb C}\setminus{\cal D}})$ elements of a
row of some fundamental matrix of a Fuchsian system with irreducible
monodromy $\chi$. Suppose, that ${\cal D}$ consists of at most three points and
the inequality
$\sum_{\cal D}(\tilde\beta_i^1+\tilde\beta_i^2)\ge 1$,
 holds. Then
$$
\left|
  \begin{array}{ccc}
    y   & \tilde y_1   & \tilde y_2   \\
    y'  & \tilde y_1'  & \tilde y_2'  \\
    y'' & \tilde y_1'' & \tilde y_2''
  \end{array}
\right|
\cdot\frac{1}
{
 \left|
   \begin{array}{cc}
     \tilde y_1  & \tilde y_2 \\
     \tilde y_1' & \tilde y_2'
   \end{array}
 \right|
}=0
$$
 is a scalar Fuchsian equation with singularities  $\cal D$ and monodromy $\chi$.
\end{Lemma}
 {\bf Proof:} The equation we study is a scalar regular equation, and hence,
it is Fuchsian equation. Suppose that zeros of the determinant give us some
additional singularities ${\cal D}_k$.
From the Fuch's relation (statement \ref{fuks}) we have
$$
\sum_{\cal D}(\tilde\beta_i^1+\tilde\beta_i^2)+
\sum_{{\cal D}_k}(\tilde\beta_i^1+\tilde\beta_i^2)\le k+1
$$
One can see that the solutions are holomorphic at points
${\cal D}_k$ and therefore the exponents at that points are nonnegative integers.
 Moreover the determinant in denominator can equal to zero at some point if and
only if there exists a linear combination of $\tilde y_1,\tilde y_2$ having zero
of at least second order at that point (the columns of the determinant are dependent).
So, for any point $b_i$ in ${\cal D}_k$ at least one of the valuations, and hence
one of the exponents is greater or equal to 2.
Then
$$
1+k\ge
\sum_{\cal D}(\tilde\beta_i^1+\tilde\beta_i^2)+
\sum_{{\cal D}_k}(\tilde\beta_i^1+\tilde\beta_i^2)
\ge 1+2k
$$
That estimate is possible only if $k=0$. That means that the scalar equation
has no additional singularities.
$\Box$
\begin{Corollary} The equation constructed with $y_1,y_2$
as we have described is the Riemann equation with the given irreducible
monodromy.
\end{Corollary}
{\bf Proof:}
From the form of $y_1,y_2$ it follows that the equation is Fuchsian.
Following the lemma it has exactly three singular points. So it is a Riemann equation.
 $\Box$
\begin{Corollary}\label{pok} The exponents of the constructed equation are
related to the exponents of the initial system as follows:
$$
\tilde\beta_j^k=\beta_j^k\qquad\mbox{при $j\ne i$ }
$$
$$
\tilde\beta_i^1=\beta_i^1+1,\,\,\,\,\tilde\beta_i^2=\beta_i^2
$$
\end{Corollary}
So the main statement of this section is proved:
\begin{Theorem}
{ Any irreducible representation of rank 2 with two generators
can be realized by the monodromy of a Riemann equation with given singular points.}
 \end{Theorem}

\subsection{Reducible representations}
Since here we consider the representations defined by two $2\times 2$,
matrices we have three cases to explore. These are diagonal, diagonalizible
and non-diagonalizible cases. The cases differ by the possibility to
diagonalize simultaneously both, one or none of the generator matrices.

\subsubsection{Diagonal monodromy}\label{diag}
In this section we consider representations that can be reduced to
$G_i=diag(\alpha_i,\beta_i)$ for all $i$ simultaneously.
For such monodromy the following lemma holds.
\begin{Lemma}\label{diagr} If the monodromy of Riemann equation is diagonal, then
at least one of $G_i$ is scalar ($\alpha_i=\beta_i$ for at least
one $i$).
\end{Lemma}
{\bf Proof:} Suppose none of the matrices $G_i$ is scalar and assume
$a_3=\infty$, that can be achieved by conformal transform of Riemann sphere.
Since there are no scalar monodromy matrix we can find two independent solutions
 and construct the global Levelt's basis:
\begin{eqnarray}
u_1(z)=(z-a_1)^{\beta_1^1}(z-a_2)^{\beta^1_2}P_1(z)\nonumber \\
  u_2(z)=(z-a_1)^{\beta_1^2}(z-a_2)^{\beta_2^2}P_2(z)\nonumber
\end{eqnarray}

Here  $\beta_i^1\ne\beta_i^2$\, and $P_1,P_2$ are polynomials of degrees $n_1,n_2$,
having no zeros at $a_1,a_2$. As the monodromy at infinity has
pairwise-distinct eigenvalues then we have
$$
\beta_3^j=-\beta_1^j-\beta_2^j-n_j
$$
That follows from the fact that $(u_1,u_2)$ is a global Levelt's basis.
Once again, we note that if the monodromy eigenvalues of diagonal monodromy
 are different in all points then any local Levelt's basis is also a global
Levelt's basis.
But from the Fuchs relation (\ref{fuks}) we have
$$
\sum_{\cal D}\left(\beta_i^1+\beta_i^2\right)=-n_1-n_2=1
$$
That contradicts to the fact that $n_1,n_2$ are non-negative. $\Box$

So, we can assume $G_3=\gamma\cdot I$. Let $y_1,y_2$ be some solutions of
a first order Fuchsian equations having monodromy
$\alpha_1,\alpha_2,\gamma$ and $\alpha_2,\beta_2,\gamma$:
\begin{eqnarray}
y_1=(z-a_1)^{\lambda_1}(z-a_2)^{\mu_1}(z-a_3)^{\nu}\nonumber\\
y_2=(z-a_1)^{\lambda_2}(z-a_2)^{\mu_2}(z-a_3)^{\nu}\nonumber
\end{eqnarray}
such that $\lambda_i+\mu_i+\nu_i=0$ (here we suppose that $\infty$
is a holomorphic point of the equation). Such a choice is always possible
because this sum is integer
\mbox{$\left(e^{2\pi i
(\lambda_i+\mu_i+\gamma)}=\alpha_i\beta_i\gamma=1\right) $},
and the monodromy define all these numbers only up to an integer.
This is also the reason for assuming \mbox{$y_1,y_2$} to be linearly independent.
 Take now \mbox{$y_1,y_2$} as the basis of solution space of some linear equation.
One can easily see that the exponents of that equation at
 \mbox{$a_1,a_2,a_3$} are \mbox{$(\lambda_1,\lambda_2),\,(\mu_1,\mu_2),\,(\nu+1,\nu) $}.
Indeed, as $y_1$ and  $y_2$ are proportional to $(z-a_3)^{\nu}$ near $z=a_3$,
 there exists their linear combination proportional to $(z-a_3)^{\nu+1}$.
Now from the lemma \ref{postrur} we obtain, that the constructed equation is a Riemann equation.
So the following result is established.
\begin{Theorem} Diagonal monodromy representation can be realised by a
Riemann equation iff one of the monodromy matrices is scalar.
\end{Theorem}
\subsubsection{Diagonalisible in a point monodromy}
By diagonalisible at a point monodromy we mean the following representation:
$$
\chi:\qquad
G_1=\left(\begin{array}{cc} g^1_{11} & 0 \\ 0 & g_{2
2}^1\end{array}\right), \,G_2=\left(\begin{array}{cc}g_{11}^2&g_{12}^2
\\ 0&g_{22}^2\end{array}\right),\,G_3=\left(\begin{array}{cc}
g_{11}^3&g_{12}^3\\0&g_{22}^3\end{array} \right)
$$
\begin{Lemma} Let
$Y(z)$ be a fundamental matrix of a fuchsian system with monodromy $\chi$.
Then the pairs \mbox{$(y_{11},y_{12})$} and \mbox{$(y_{21},y_{22})$}
can not be linearly dependent simultaneously.
\end{Lemma}
{\bf Proof:}
Suppose $k\cdot y_{11}=y_{12}$\,,$\,y_{21}=l\cdot y_{22}$.
We have $kl\ne 1$ because of the linear independence of columns of $Y(z)$.
Then in the basis
$$
\widetilde Y=Y
\left(
   \begin{array}{cc}
      \frac{1}{1-kl} & \frac{-k}{1-kl} \\
      \frac{-l}{1-kl} & \frac{1}{1-kl}
   \end{array}
\right)
=
\left(
    \begin{array}{cc}
      \tilde y_1 & 0 \\
      0 &\tilde y_4
    \end{array}
\right)
$$
the fundamental matrix is diagonal, and hence the monodromy (global) is diagonal.
But $\chi$ is not diagonal. So a simultaneous linear dependence is impossible. $\Box$

Now let $Y(z)$ be the fundamental matrix constructed in a Levelt's basis at
the point with diagonal monodromy. We can that this point is $a_1$.
Consider
$$
\widetilde Y=
\left(
    \begin{array}{cc}
       0&1\\
       1&0
    \end{array}
\right)
U_1^{-1}(a_1)Y
$$
For $\tilde{Y}(z)$ we have
$$
U_1(a_1)=\left(\begin{array}{cc} 0 & 1 \\1 & 0 \end{array}\right)
$$
and either \mbox{$(y_{11},y_{12})$}, or \mbox{$(y_{21},y_{22})$}
are linearly independent. Considering an equation having independent pair as
the basis of solution space, we see that this pair is either Levelt's basis of
the equation at $a_i$, or its transposition ($G_i$ is diagonal).
So the conditions of lemma (\ref{postrur}) hold, and the equation obtained
is a Riemann equation.
\begin{Theorem}
Any reducible indecomposable representation diagonalasible at least at
one point admits realisation by the monodromy of a Riemann equation with given
singular points.
\end{Theorem}
\subsubsection{Non-diagonalisible monodromy}
In this section we deal with the monodromy of the following type
$$
\chi:\qquad G_1=\left(\begin{array}{cc}\lambda_1 &
1\\0&\lambda_1\end{array}\right),\,G_2=\left(\begin{array}{cc}\lambda_2&
c\\0&\lambda_2\end{array}\right),\,G_3=\left(\begin{array}{cc}\lambda_3&d
\\0&\lambda_3\end{array}\right)
$$
Consider a basis $(y_1,y_2)$ of the solution space of a Riemann equation
with the monodromy of a given form. The basis $(y_1,y_2)$ is a global
Levelt's basis. Indeed, the only invariant subspace of the operator having the
Jordan block matrix is the line spanned by the first vector of Jordan basis.
Hence all exponents of the equation can be defined in that basis.
One can note, that \mbox{$\beta_i^1$} is a logarithmic residue of $y_1$ at
a point $a_i$. The sum of all residues of $y_1$ over ${\mathbb C}P^1$ equals
to zero and consists of the sum of residues over ${\cal D}$ and non-negative sum
over zeros of $y_1$. Therefore \mbox{$\sum_{\cal D}\beta_i^1\le 0$}.
As the basis in consideration is the Levelt's basis then we have
\mbox{$\Re\beta_i^1\ge\Re\beta_i^2$}. And so
$
\sum_{\cal D}\left(\beta_i^1+\beta_i^2\right)\le 2 \sum_{\cal D}\beta_i^1\le0
$
 in contradiction to Fuchs relation (\ref{fuks}). We see that the assumption of
the existing of the basis $y_1,y_2$ was false.
\begin{Theorem} Any monodromy realisible by Riemann equation is
diagonalizible at least at one point.
\end{Theorem}

\section{$SL(2,{\mathbb C})$ monodromy}
It is known, (\cite{gol}) that any Riemann equation is completly defined by
its divisor ${\cal D}$ and the set of exponents $\beta_i^j$. More precisely, for
\mbox{$\infty\not\in {\cal D}=\{a_1,a_2,a_3\}$}, we have:
\begin{equation}
y''+
\left(
\sum_{\cal D}\frac{1-\beta_i^1-\beta_i^2}{z-a_i}
\right)
y'+
\left(
\sum_{\cal D}\frac{\beta_i^1\beta_i^2}{z-a_i}\cdot
\frac{\prod_{j\ne i}(a_i-a_j)}{\prod_{\cal D}(z-a_i)}
\right)
y=0 \label{aa}
\end{equation}
for a Riemann equation having the exponents $\beta_i^j$.
If \mbox{$\infty\in{\cal D}$} then assuming $a_3=\infty$ we get
\begin{eqnarray}
y''+\left(\sum_{i=1,2}\frac{1-\beta_i^1-\beta_i^2}{z-a_i}
\right)y'+\qquad\qquad\qquad\qquad\qquad\qquad\qquad\qquad\qquad\qquad
\qquad \nonumber \\ +\left(\frac{\beta_1^1\beta_1^2(a_1-a_2)}{z-a_1}
+\frac{\beta_2^1\beta_2^2(a_2-a_1)}{z-a_2}+\beta_3^1\beta_3^2\right)
\frac{y}{(z-a_1)(z-a_2)}=0\label{bb}
\end{eqnarray}
\begin{Def} The Rieman-Sturm-Liouville equation is a Riemann equation
 having no term with $y'$:
$$
 y''+q(z)y=0
$$
\end{Def}
From (\ref{aa}),(\ref{bb}) and the Fuchs relation on $\beta_i^j$ it follows
 that Riemann-Sturm-Liouville equation has $\infty$  as one of its singular points.
Let us discuss the monodromy of such an equation.
We have for  $i=1,2$
$$
1=\beta_i^1+\beta_i^2=(\rho_i^1+\rho_i^2)+(\varphi_i^1+\varphi_i^2)\Rightarrow
 (\rho_i^1+\rho_i^2)\in{\mathbb Z}\Rightarrow 1=e^{2\pi
i(\rho_i^1+\rho_i^2)}=\det G_i
$$
So the monodromy of RSL equation lies in $SL\left(2,{\mathbb C}\right)$.
There naturally arises the question whether all possible $SL$ monodromies of
Riemann equation admit the realization by RSL equation too.
\begin{Sta}
\label{realis} The representation \mbox{$\chi:\pi_1(\overline{\mathbb{C}}
\setminus{\cal D})\to SL(2,C)$} admits the realization by the RSL equation monodromy
 iff there exists  a Riemann equation with the same monodromy such that
$\beta_i^1+\beta_i^2=1$ for two of its singular points.
\end{Sta}
{\bf Proof:} Suppose that at the points $a_1,a_2$ the sum of exponents
equals to one. By a conformal map of Riemann sphere one can move $a_3$ to
infinity. Following (\ref{bb}) the resulting equation will be RSL equation.
The converse is trivial. $\Box$

It remains to find out which of $SL$ monodromies can be Realized by
described Riemann equation. As in Section 1 we will distinguish between
different representation classes.

\subsection{Diagonal monodromy}
First we realize the degenerate case. Set $G_i=I$ for all $i$.
 To that end consider
\begin{eqnarray} y_1=(z-a_1)^2(z-a_2)^{-1}(z-a_3)^{-1}\nonumber\\
                 y_2=(z-a_1)^{-1}(z-a_2)^2(z-a_3)^{-1}\nonumber
\end{eqnarray}
As it was proved in \ref{diag}, that is the basis of some Riemann equation's solution space.
From the corollary \ref{pok} we have
$\beta_i^1+\beta_i^2=\varphi_i^1+\varphi_i^2=1\,$ for $i=1,2$. From the statement 2.1.
we obtain that this Riemann equation is in fact RSL equation. Now,
if not all $G_i$ are unit matrices, then one can check that
 $\rho_1^j+\rho_2^j+\rho_3^j=1\,,\,j=1,2$ holds. Here $\rho_i^j$ are
the normalized logarithms of the monodromy. Assuming according to lemma
\ref{diagr}, that $G_3=\pm I$ we choose $y_1,y_2$ in the following way:
\begin{eqnarray}
y_1=(z-a_1)^{k+\rho_1^1}\,\cdot\,(z-a_2)^{-k+\rho_2^1}\cdot(z-a_3)^{-1+\rho_3}\nonumber\\
y_2=(z-a_1)^{-k+\rho_1^2}\cdot(z-a_2)^{k+\rho_2^2}\,\cdot\,(z-a_3)^{-1+\rho_3}\nonumber
\end{eqnarray}
As $G_i\ne I$ at least in two points then  we have
 $\beta_i^1+\beta_i^2=\rho_i^1+\rho_i^2=1$. From that follows it that
the Riemann equation we constructed satisfy the conditions of statement \ref{pok}.
Hence any diagonal monodromy with at least one of $G_i$'s being scalar can be
realized by Riemann-Sturm-Liuoville equation.

 \subsection{Other representations}
In this section we will need some additional facts from the Fuchs systems
theory. Consider a fuchsian system having no singular point at infinity.
Such a system has the coefficient matrix of the following form:
$$
B(z)=\sum_{\cal D}\frac{B_i}{z-a_i},\,\, \sum B_i=0
$$
\begin{Lemma} Let $Y$ be a fundamental matrix of fuchsian system
\mbox{$\dot Y=B(z)Y$} constucted in a Levelt's basis at $a_i$
$$
Y=U_i(z)(z-a_i)^{A_i}(z-a_i)^{E_i}
$$
Denote $L_i=\lim_{z\to{a_i}}(z-a_i)^{A_i}E_i(z-a_i)^{-A_i}+A_i$.
 Then
$$
B_i=U_i(a_i)L_iU_i^{-1}(a_i)
$$
\end{Lemma}
{\bf Proof:} In a neighborhood of $a_i$ we have:
$$
B(z)=\frac {dY}{dz} Y^{-1}=\left(
\frac{dU_i(z)}{dz}U_i^{-1}(z)+
\frac{U_i(z)
         \left( A_i+(z-a_i)^{A_i}E_i(z-a_i)^{-A_i}\right)
     U_i^{-1}(z)}
   {z-a_i}\right)
$$
That immediately proves the lemma.  $\Box$

\noindent {\bf Замечание:} If $A_i$ is a scalar matrix, then $L_i=A_i+E_i$

The main result this section is based on is the following lemma.
\begin{Lemma} (\cite{Bol}) For any $\chi:\pi_1(\overline{\mathbb C}\setminus {\cal
D}\to GL(2,{\mathbb C}))$ there exists a Fuchsian system with a given monodromy such
that $\varphi_i^1=\varphi_i^2=0$ in at least in two points.
\end{Lemma}
{\bf Proof:} \label{obnul} Suppose there exist two points with coincident
valuations. For example, \mbox{$\varphi_1^1=\varphi_1^2=k$} and
\mbox{$\varphi_2^1=\varphi_2^2=l$} holds for points $a_1,a_2$.
In that case we turn from the solution space ${\cal X}$
to
$$
{\cal X'}=
\left(\frac{(z-a_3)^{k+l}}{(z-a_1)^k(z-a_2)^l}\right)^I{\cal
X}=
\left(
\begin{array}{cc}
\frac{(z-a_3)^{k+l}} {(z-a_1)^k(z-a_2)^l} & 0 \\
 0 & \frac{(z-a_3)^{k+l}}{(z-a_1)^k(z-a_2)^l}
\end{array}
\right)
{\cal X}
$$
The latter is evidently a solution space of some fuchsian system with the
same monodromy and required valuations.

Now suppose that in at least two points (suppose again these are $a_1,a_2$) we have
\mbox{$\varphi_i^1-\varphi_i^2>0$}. The following procedure will be applied.
Turn at first to
$$
{\cal X}'=U_1^{-1}(a_1){\cal X}
$$
and then to
$$
{\cal X}''=\Gamma_1(z){\cal X}=
\left(
\begin{array}{cc}
\frac{z-a_2}{z-a_1}&0\\
0&1
\end{array}
\right)
{\cal X}'
$$
From the form of $U_1''(z)$ it follows that
${\varphi_1^1}''=\varphi_1^1-1$,\,\,\,${\varphi_1^2}''=\varphi_1^2$.
Let
$$
U_2''(z)=
\left(
\begin{array}{cc}
 (z-a_2)(a_0+...)& (z-a_2)(b_0+...)\\
c_0+...& d_0+...
\end{array}
\right)
$$
and assume $a_0\ne 0$. Then taking
$$
\Gamma_2(z)=
\left(
\begin{array}{cc}
 1 & 0 \\
 -\frac{c_0a_0^{-1}}{z-a_2}&1
\end{array}
\right)
$$
we obtain
$$
U_2'''(z)=\Gamma_2(z)U_2''(z)=\left(\begin{array}{ll}(z-a_2)(a_0+...)&
(z-a_2)(b_0+..)\\ (z-a_2)(c_1-c_0a_1a_0^{-1}+...)& (a_0^{-1}\det
U_2'(a_2)+...)\end{array}\right)
$$
Obviously here, $\varphi_1^{1'''} = \varphi_1^{1} - 1$,
$\varphi_2^{1'''}=\varphi_2^1+1$. From the invariance of the sum of exponents
after all these transformations applied it follows that the resulting system
is also fuchsian. As $\Gamma_i(z)$ is meromorphic and holomorphically invertible
over $\overline{{\mathbb C}}\setminus{\cal D}$  it follows that the monodromy
is also preserved.

Now if $a_0=0$, then according to statement \ref{nevyr} we have $b_0\ne 0$.
In a similar way taking
$$
\Gamma_3(z)=
\left(
\begin{array}{cc}
1& 0 \\
-\frac{d_0b_0^{-1}}{z-a_2}&1
\end{array}
\right)
$$
in  a similar way we obtain  $\varphi_1^{1'''}=\varphi_1^1-1$\,,\,$
\varphi_2^{1'''}=\varphi_2^2+1$. All other valuations do not change.

Finally, after described procedure, the difference
 between $\varphi_i^1$ and $\varphi_i^2$ decreases at least at one point.
Applying this process until one of the differences becomes zero,
and realizing if necessary the same procedure for other residuary points
we get the already studied case where $\varphi_i^1=\varphi_i^2$, at least
at two points. The lemma is proved. $\Box$

{\bf Замечание} The gauge transformations $\Gamma_2(z)$, $\Gamma_3(z)$
we applied are special cases of more general $A$- and $B$-transformations
constructed in \cite{Bol}.

{\bf Corollary:} If $\chi$ is a non-diagonal $SL$-representation,
then exists a fuchsian system with a given monodromy representation and
prescribed singularities ${\cal D}$ such that
\mbox{$\Sigma_i=\beta_i^1+\beta_i^2=1$} in at least two points.

{\bf Proof:} Consider a system constructed in Lemma \ref{obnul}.
For $\varepsilon_i=\rho_i^1+\rho_i^2$ the following three cases
(up to permutation) have to be considered:
\begin{enumerate}
\item $\Sigma_1=\Sigma_2=1$,\,\,$\Sigma_3=-2$ --
this system is in the corollary conditions.
\item
$\Sigma_1=1$,\,$\Sigma_2=0$,\,$\Sigma_3=-1$ Turning to ${\cal
X'}=\Gamma(z){\cal X}$ where
$\Gamma(z)=\left(\frac{z-a_3}{z-a_2}\right)^I$, we get
\mbox{$\Sigma_1'=\Sigma_3'=1$}
 \item $\Sigma_1=\Sigma_2=\Sigma_3=0$ (that is
$\varepsilon_1=\varepsilon_2=\varepsilon_3=0$). We apply the following gauge:
$$
{\cal X'}=
\left(
  \begin{array}{cc}
    0&1\\
    1&0
  \end{array}
\right)
U_1^{-1}(a_i){\cal X}
$$
Here we assumed $U_1(z)$ to be holomorphicaly invertible part of the Levelt's
decomposition of the Levelt's basis (at $a_1$) matrix $Y(z)$. Assume
the upper-left element $u^i_{11}(a_i)$ of at least one of the matrices
$U_2'(a_2),\,U_3'(a_3)$ to be non-zero. In this case we can apply the
procedure of the lemma \ref{obnul} in order to increase one of the
valuations, for example $\varphi_2^1$, and decrease $\varphi^1_2$
getting $\Sigma_1'=-1$, $\Sigma_2'=1$, $\Sigma_3'=0$. Thus the problem is
reduced to the already explored sutuation.

If $u_{11}^i(a_i)=0$ for any $i$ we have
\begin{eqnarray} U_i(a_i) & =
&
\left(
\begin{array}{cc}
0&b\\
c&d
\end{array}
\right)
\nonumber\\
B_i\,\quad & =
& U_i(a_i)L_iU_i^{-1}(a_i)=
\left(
\begin{array}{rr}
0&b\\
 c&d
\end{array}
\right)
\left(
\begin{array}{cc}
0&\lambda\\
0&0\end{array}
\right)
\left(
\begin{array}{rc}
-\frac{d}{bc}&\frac{1}{c}\\
\frac{1}{b}&0
\end{array}
\right)
=
\left(
\begin{array}{cc}
0&0\\
\frac{\lambda c}{b}&0
\end{array}
\right)
\nonumber
\end{eqnarray}
and $B(z)=\sum\frac{B_i}{z-a_i}$ is of the following form
$$
B(z)=\left(\begin{array}{cc} 0&0\\f(z)&0\end{array}\right)
$$
Such a system is reducible, and its monodromy is also reducible, because of
the existence of solution of the form
$y=\left(\begin{array}{c}0\\c\end{array}\right)$.
That solution is invariant under analytical continuation along any loop.
So one can pose $g^i_{21}=0$ for all $i$ in the basis having $y$ as the first
 basis vector. As it is proved in the section devoted to $GL(2,{\mathbb C})$
monodromy one of the monodromy matrices $G_i$ can be assumed being diagonal.
And now we can decrease one of the valuations at $a_j$ increasing one of
them at $a_i$. Since $G_i$ is diagonal, then the basis we obtained is either Levelt's
basis, or its transposition. In that basis we have $\Sigma_1=0$, $\Sigma_2=-1$,
$\Sigma_3=1$ -- the case which already have been considered.
 $\Box$

\end{enumerate}

Now we are ready to prove the existence of Riemann equation  with the monodromy
we are interested in and given exponents sum.

Suppose $\chi$ is irreducible. There exists such a system, having $\chi$
 as the monodromy, that $\Sigma_1=\Sigma_2=1$. Proceeding to the equation
by increasing the valuations at $a_3$ we get according to corollary \ref{pok}
the desired Riemann equation.

If $\chi$ is reducible then one of $G_i$ can be considered diagonal, and one
have to increase the valuations exactly at that point. If that is $G_3$ then
in the way similar to irreducible case we get the Riemann equation. If that
is $G_2$ then we increase $\varphi_3^1$ decreasing one of $\varphi_2^j$. More
precisely we decrease $\varphi_3^2$ increasing one of $\varphi_2^j$ and
after that we diminish by 2 the exponents sum at $a_2$ and enlarge it by 2
at $a_3$ as in the first part of the lemma. After that, as the matrix $G_2$
is diagonal, the image of initial Levelt's basis is either the Levelt's basis
of the resulting equation or its transposition. In any case by the corollary
 \ref{pok}  the sum of exponents both at $a_1$ and  $a_2$ is equal to one.
\begin{Theorem} Any $SL(2,{\mathbb C})$ representation of the monodromy of Riemann
equation can be realized by the monodromy of Riemann-Sturm-Liouville
equation as well.
\end{Theorem}

One can note that in an important special case of hypergeometric equation the
condition for monodromy to be in $SL(2,{\mathbb C})$ is especially
simple. The hypergeometric equation
$$
 u''+\frac{\gamma-(1+\alpha+\beta)z}{z(z-1)}u'-\frac{\alpha\beta }{z(z-1)}u=0
$$
has the monodromy in $SL(2,{\mathbb C})$ iff  $\{\gamma,\alpha+\beta\}\in
{\mathbb
Z}$.

\section{$SL(2, {\mathbb Z})$ monodromy of hypergeometric equation}

It is evident that any {$SL(2, {\mathbb Z})$ monodromy realized by Riemann
equation and having eigenvalues equal to one for at least two monodromy matrices
can be realized by hypergeometric equation. One obtain it from that Riemann
equation by a conformal map of the Riemann sphere.
There arises natural question, on given hypergeometric equation,
to define whether its monodromy is in $SL(2, {\mathbb Z})$.
The obvious necessary condition is the belonging of the monodromy to
$SL(2, {\mathbb C})$ . In terms of the coefficients of the equation it is
$\{\alpha+\beta, \gamma\}\in{\mathbb Z}$.
Let us elucidate whether these conditions are sufficient.
One can easily check that for reducible monodromy representations this
conditions are indeed sufficient. In fact any reducible monodromy of
a hypergeometric equation can be reduced to
$$
G_0=
\left(
\begin{array}{cc}
1&c\\
0&1
\end{array}
\right)
,\quad
G_1=
\left(
\begin{array}{cc}
1&0\\
0&1
\end{array}
\right)
,\quad
G_{\infty}=G_1^{-1}G_0^{-1}
$$
After conjugation by $S_0=diag(c^{-1},1)$ we have
$$
G_0=
\left(
\begin{array}{cc}
1&1\\
0&1
\end{array}
\right)
,\quad
G_1=
\left(
\begin{array}{cc}
 1&0\\
0&1
\end{array}
\right)
,\quad
G_{\infty}=G_1^{-1}G_0^{-1}
$$

In irreducible case the further examination is required. Any irreducible
$SL$ representation of a hypergeometric equation can be reduced to
$$
G_0=
\left(
\begin{array}{cc}
1&0\\
d&1
\end{array}
\right)
,\quad
G_1=\left(
\begin{array}{cc}
1&c\\
0&1
\end{array}
\right)
,\quad
G_{\infty}=G_1^{-1}
G_0^{-1},\qquad c,d\in\mathbb C
$$
Conjugation by  $S_1=diag(d,1)$ gives
$$
G_0=
\left(
\begin{array}{cc}
1&0\\
1&1
\end{array}
\right)
,\quad
G_1=
\left(
\begin{array}{cc}
1&cd\\
0&1
\end{array}
\right)
,\quad
cd\in\mathbb C
$$
Now after conjugation by an arbitrary $S_2 \in SL(2,\mathbb C)$ we obtain
$$
G_0=
\left(
\begin{array}{cc}
1+\mu\tau&-\mu^2\\
 \tau^2&1-\mu\tau
\end{array}\right)
,\quad
G_1=
\left(
\begin{array}{cc}
 1-\lambda\nu b& \lambda^2b\\
-\nu^2b&1+\lambda\nu b
\end{array}
\right)
,\quad b=cd\in\mathbb C
$$
where
$$
S_2=
\left(
\begin{array}{cc}
\lambda&\mu\\
 \nu&\tau
\end{array}
\right)
$$
One can see that the criterion of possibility to reduce the monodromy to
$SL(2,{\mathbb Z})$ is the existence of the set $\{\lambda,\mu,\nu,\tau\}$
satisfying the following relations:
$$
\left\{
\begin{array}{l}
\mu^2,\tau^2,\mu\tau\in{\mathbb Z}\\
\lambda^2b,\nu^2b,\lambda\nu b\in{\mathbb Z}\\
\lambda\tau-\mu\nu=1
\end{array}
\right.
$$
If such a set exists, then
$$
\lambda^2\tau^2b,\mu^2\nu^2b,\lambda\mu\nu\tau b\in{\mathbb Z}\Rightarrow(\lambda^2\tau^2+\mu^2\nu^2-2\lambda\mu\nu\tau )b=
(\lambda\tau-\mu\nu)^2b=b\in{\mathbb Z}
$$
One can note that $b=tr G_{\infty}-2=e^{2\pi i\alpha}+e^{2\pi i\beta}-2$.
Therefore we get the following statement
\begin{qwe} The hypergeometric equation monodromy belongs
 to $SL(2,{\mathbb Z})$ if and only if $\{\gamma,\alpha+\beta,
e^{2\pi i\alpha}+e^{2\pi i\beta}\}\in{\mathbb Z}$ holds.
\end{qwe}

This conditions can be examinated in greater details. We have
$$
e^{2\pi i\alpha}+e^{2\pi i\beta}=e^{2\pi i\alpha}(1+e^{2\pi i(\alpha+\beta-2\alpha)})=
e^{2\pi i\alpha}+e^{-2\pi i\alpha}=e^{2\pi i\beta}+e^{-2\pi
i\beta}=k,\phantom{aa}k\in\mathbb{Z}
$$
And thus for $x=e^{2\pi i\alpha}$ we get $x+x^{-1}=k$. Together
with $(\alpha+\beta)\in\mathbb{Z}$ that gives
$e^{2\pi i\alpha}=\frac{k+\sqrt{k^2-4}}2,\phantom{a}e^{2\pi i\beta}=\frac{k-\sqrt{k^2-4}}2$

Finally, the monodromy of hypergeometric equation lies in
$SL(2,\mathbb{Z})$ for the following parameter's values:
$$
\alpha=\frac 1{2\pi i}\ln\frac {k\pm\sqrt{k^2-4}}2,\phantom{a}\beta=-\alpha=
\frac 1{2\pi i}\ln\frac{k\mp\sqrt{k^2-4}}2,\phantom{a}\gamma=l,\phantom{a}
k,l\in\mathbb{Z}
$$

The corresponding equation appears to have following form
$$
 u''+\frac{l-z}{z(z-1)}u'-\frac{\frac 1{(2\pi)^2}\left(\ln\frac{k\pm\sqrt{k^2-4}}{2}\right)^2 }{z(z-1)}u=0
$$
with any integer $k,l$.


\newpage
\begin{thebibliography}{99}

\bibitem{b7} Bolibruch~A.A {\it On sufficient conditions of positive
solvability of Riemann-hilbert problem} // Math. Notes,1992, vol 51,2.

\bibitem{Bol} Bolibruch~A.A {\it The 21st Hilbert problem for Linear
Fuchsian Systems.} Moscow, Nauka, 1994

\bibitem{dek} Dekkers W. {\it The matrix of a connection having regular
singularities on a vector bundle of rank 2 on
$P^1(C)$}//Lect.Notes.Math.1979, vol 712.

\bibitem{for} Forster~O. {\it Riemann surfaces.} Moscow, Mir,1980.

\bibitem{Gan} Gantmacher~F.R. {\it Theory of Matrices}, vol 2, New York,
Chelsea, 1959

\bibitem{gol} Golubev~V.V. {\it Lections on analytical theory of
differential equations.}  Moscow, Gostechteorizdat,1950.

\bibitem{yap} Iwasaky K., Kimura H., Shimomura Sh., Masaaky Y.
{\it From Gauss to Painlev\'e.}
 Vieweg verlag, 1991, Aspects of Math. E; vol 16.

\bibitem{kim} Kimura,~T. and Shima,~K. A note on the monodromy of the
hypergeometric differential equation, Japan.~J.~Math., 17 (1991), 137-163

\bibitem{kr} Krylov~B.L. {\it Explicite solution of riemann problem for
Gauss system.} // Kazan, Tr. Kaz. Av. inst.,1956,vol 31.

\bibitem{lev} Levelt A.H.M. {\it Hy\-per\-ge\-o\-met\-ric func\-ti\-ons.}
// Proc. Knkl. netherl. Acad.wetensch., Ser A., 1961 vol 64.

\bibitem{zog} Zograf~P.G. and Tachtadjan~L.A
{\it On the Liouville equation, accessory parameters and the geometry
of the Teichm\"uller space for the Riemann surfaces of genus 0.} //Math.
Sb.,1987, vol 132,2.
\end{thebibliography}
\end{document}